\documentstyle{article}
\begin{document}
\title{\Large{Lie symmetries and conservation laws of the Hirota-Ramani equation }}
\author{Mehdi Nadjafikhah\thanks{School of Mathematics, Iran University of Science and Technology, Narmak, Tehran 1684613114, Iran. e-mail:
m\_nadjafikhah@iust.ac.ir} \and Vahid
Shirvani-Sh.\thanks{Department of Mathematics, Islamic Azad
University, Karaj Branch, Karaj 31485-313, Iran. e-mail:
v.shirvani@kiau.ac.ir } }
\date{}
\maketitle
\begin{abstract}
In this paper, Lie symmetry method is performed for the
Hirota-Ramani (H-R) equation. We will find The symmetry group and
optimal systems of Lie subalgebras. Furthermore, preliminary
classification of its group invariant solutions, symmetry
reduction and nonclassical symmeries are investigated. Finally
the conservation laws of the H-R equation are presented.
\end{abstract}


\input{amssym}
\def\be{\begin{eqnarray}}
\def\ee{\end{eqnarray}}
\def\di{\displaystyle}
\def\rank{{\bf rank}}
\section{Introduction}
In the present paper, we study the following equation
\be u_t-u_{x^2t}+au_{x}(1-u_{t})=0,\label{eq:1.0}\ee
where $a\neq0$ and $u(x,t)$ is the amplitude of relevant wave
mode. This equation was introduced by Hirota and Ramani in
\cite{[8]}. Jie Ji obtained some travelling soliton solutions of
this equation by using Exp-function method \cite{[15]}. This
equation is completely integrable by the inverse scattering
method. Eq. (\ref{eq:1.0}) is studied in \cite{[8],[15],[10]}
where new kind of solutions were obtained. Hirota-Ramani equation
is widely used in various branches of physics, such as plasma
physics, fluid physics, quantum field theory. It also describes a
variety of wave phenomena in plasma and solid state \cite{[8]}.

The theory of Lie symmetry groups of differential equations was
developed by Sophus Lie \cite{[2]}, which was called classical
Lie method. Nowadays, application of Lie transformations group
theory for constructing the solutions of nonlinear partial
differential equations (PDEs) can be regarded as one of the most
active fields of research in the theory of nonlinear PDEs and
applications. Such Lie groups are invertible point transformations
of both the dependent and independent variables of the
differential equations. The symmetry group methods provide an
ultimate arsenal for analysis of differential equations and is of
great importance to understand and to construct solutions of
differential equations. Several applications of Lie groups in the
theory of differential equations were discussed in the
literature, the most important ones are: reduction of order of
ordinary differential equations, construction of invariant
solutions, mapping solutions to other solutions and the detection
of linearizing transformations for many other applications of Lie
symmetries see \cite{[3],[4],[7]}.

The fact that symmetry reductions for many PDEs are unobtainable
by applying the classical symmetry method, motivated the creation
of several generalizations of the classical Lie group method for
symmetry reductions. The nonclassical symmetry method of
reduction was devised originally by Bluman and Cole in 1969
\cite{[11]}, to find new exact solutions of the heat equation. The
description of the method is presented in \cite{[9],[6]}. Many
authors have used the nonclassical method to solve PDEs. In
\cite{[12]} Clarkson and Mansfield have proposed an algorithm for
calculating the determining equations associated to the
nonclassical method. A new procedure for finding nonclassical
symmetries has been proposed by B\^{i}l\v{a} and Niesen in
\cite{[5]}.

Many PDEs in the applied sciences and engineering are continuity
equations which express conservation of mass, momentum, energy,
or electric charge. Such equations occur in, e.g., fluid
mechanics, particle and quantum physics, plasma physics,
elasticity, gas dynamics, electromagnetism,
magneto-hydro-dynamics, nonlinear optics, etc. In the study of
PDEs, conservation laws are important for investigating
integrability and linearization mappings and for establishing
existence and uniqueness of solutions. They are also used in the
analysis of stability and global behavior of solutions
\cite{[16],[17],[18],[19]}.

This work is organized as follows. In section 2 we recall some
results needed to construct Lie point symmetries of a given system
of differential equations. In section 3, we give the general form
of a infinitesimal generator admitted by Eq. (\ref{eq:1.0}) and
find transformed solutions. In Section 4, we construct the optimal
system of one-dimensional subalgebras. Lie invariants, similarity
reduced equations and differential invariants corresponding to
the infinitesimal symmetries of Eq. (\ref{eq:1.0}) are obtained
in section 5 and 6. Section 7, is devoted to the nonclassical
symmetries of the H-R model, symmetries generated when a
supplementary condition, the invariance surface condition, is
imposed. Finally in last section, the conservation laws of the
Eq. (\ref{eq:1.0}) are obtained.
\section{Method of Lie Symmetries}
In this section, we recall the general procedure for determining
symmetries for any system of partial differential equations see
\cite{[3],[13],[4],[7]}. To begin, let us consider the general
case of a nonlinear system $E$ of partial differential equations
of order $n$ in $p$ independent and $q$ dependent variables is
given as a system of equations
\be \Delta_\nu(x,u^{(n)})=0,\;\;\;\;\; \nu=1,\cdots,l,
\label{eq:2.1} \ee
involving $x = (x^1,\cdots, x^p)$, $u = (u^1,\cdots,u^q)$ and the
derivatives of $u$ with respect to $x$ up to $n$, where $u^{(n)}$
represents all the derivatives of $u$ of all orders from $0$ to
$n$. We consider a one-parameter Lie group of infinitesimal
transformations acting on the independent and dependent variables
of the system (\ref{eq:2.1})
\be \tilde{x}^i &=& x^i+s \xi^i(x,u)+O(s^2), \hspace{1cm}
i=1\cdots,p,\nonumber \\[-2mm] \label{eq:2.2}\\[-2mm] \tilde{u}^j &=& u^j+s
\varphi^j(x,u)+O(s^2), \hspace{9mm} j=1\cdots,q, \nonumber
 \ee
where $s$ is the parameter of the transformation and $\xi^i$,
$\eta^j$ are the infinitesimals of the transformations for the
independent and dependent variables, respectively. The
infinitesimal generator ${\mathbf v}$ associated with the above
group of transformations can be written as
\be  {\mathbf v} = \sum_{i=1}^p\xi^i(x,u)\partial_{x^i} +
\sum_{\alpha=1}^q\varphi^\alpha(x,u)\partial_{u^\alpha}.
\label{eq:2.4} \ee
A symmetry of a differential equation is a transformation which
maps solutions of the equation to other solutions. The invariance
of the system (\ref{eq:2.1}) under the infinitesimal
transformations leads to the invariance conditions (Theorem 2.36
of \cite{[3]})
\be \textrm{Pr}^{(n)}{\mathbf
v}\big[\Delta_\nu(x,u^{(n)})\big]=0,\;\;\;\;\;
\nu=1,\cdots,l,\;\;\;\;\mbox{whenever}\;\;\;\;\;\Delta_\nu(x,u^{(n)})
=0, \label{eq:2.5} \ee
where $\textrm{Pr}^{(n)}$ is called the $n^{th}$ order
prolongation of the infinitesimal generator given by
\be \textrm{Pr}^{(n)}{\mathbf v}= {\mathbf
v}+\sum^q_{\alpha=1}\sum_J
\varphi^J_\alpha(x,u^{(n)})\partial_{u^\alpha_J},\label{eq:2.6}
\ee
where $J=(j_1,\cdots,j_k)$, $1\leq j_k\leq p$, $1\leq k\leq n$ and
the sum is over all $J$'s of order $0<\# J\leq n$. If $\#J=k$, the
coefficient $\varphi_J^\alpha$ of $\partial_{u_J^\alpha}$ will
only depend on $k$-th and lower order derivatives of $u$, and
\be \varphi_\alpha^J(x,u^{(n)})=D_J(\varphi_\alpha - \sum_{i=1}^p
\xi^iu_i^\alpha) + \sum_{i=1}^p\xi^iu^\alpha_{J,i}, \label{eq:2.7}
\ee
where $u_i^\alpha:=\partial u^\alpha/\partial x^i$ and
$u_{J,i}^\alpha := \partial u_J^\alpha/\partial x^i$.

\medskip One of the most important properties of these
infinitesimal symmetries is that they form a Lie algebra under the
usual Lie bracket.
\section{Lie symmetries of the H-R equation }
We consider the one parameter Lie group of infinitesimal
transformations on $(x^1=x,x^2=t,u^1=u)$,
\be \tilde{x} &=& x+s\xi(x,t,u)+O(s^2),\nonumber\\
\tilde{t} &=& x+s\eta(x,t,u)+O(s^2),\label{eq:3.1}\\
\tilde{u} &=& x+s\varphi(x,t,u)+O(s^2),\nonumber \ee
where $s$ is the group parameter and $\xi^1=\xi$, $\xi^2=\eta$ and
$\varphi^1=\varphi$ are the infinitesimals of the transformations
for the independent and dependent variables, respectively. The
associated vector field is of the form:
\be {\mathbf
v}=\xi(x,t,u)\partial_x+\eta(x,t,u)\partial_t+\varphi(x,t,u)\partial_u.
\label{eq:3.2}\ee
and, by (\ref{eq:2.6}) its third prolongation is
\be \textrm{Pr}^{(3)}{\mathbf v} &=& {\mathbf v}+
\varphi^x\,\partial_{u_x}+\varphi^t\,\partial_{u_t}+\varphi^{x^2}\,\partial_{u_{x^2}}
+\varphi^{xt}\,\partial_{u_{xt}}+\varphi^{t^2}\,\partial_{u_{t^2}} \nonumber\\
&&
+\varphi^{x^3}\,\partial_{u_{x^3}}+\varphi^{x^2t}\,\partial_{u_{x^2t}}+\varphi^{xt^2}\,\partial_{u_{xt^2}}+\varphi^{t^3}\,\partial_{u_{t^3}}.
\label{eq:3.2-1}\ee
where, for instance by (\ref{eq:2.7}) we have
\be
\varphi^x&=&D_x(\varphi-\xi\,u_x-\eta\,u_t)+\xi\,u_{x^2}+\eta\,u_{xt},\nonumber\\
\varphi^t&=&D_t(\varphi-\xi\,u_x-\eta\,u_{t})+\xi\,u_{xt}+\eta\,u_{t^2},\nonumber\\
&& \vdots \label{eq:3.3} \\
\varphi^{t^3}&=&D^3_{t}(\varphi-\xi\,u_x-\eta\,u_t)+\xi\,u_{xt^3}+\eta\,u_{t^4},\nonumber
\ee
where $D_x$ and $D_t$ are the total derivatives with respect to
$x$ and $t$ respectively.
By (\ref{eq:2.5}) the vector field ${\mathbf v}$ generates a one
parameter symmetry group of Eq. (\ref{eq:1.0}) if and only if
\be \left. \begin{array}{l} \di \textrm{Pr}^{(3)}{\mathbf
v}[u_t-u_{x^2t}+au_{x}(1-u_{t})]=0,\\[5mm]
\di \mbox{whenever} \hspace{.5cm} u_t-u_{x^2t}+au_{x}(1-u_{t})=0.
\end{array}\right. \label{eq:3.3-1} \ee
The condition (\ref{eq:3.3-1}) is equivalent to
\be \left. \begin{array}{l} \di (1-au_{x})\varphi^t+a(1-u_{t})\varphi^x-\varphi^{x^2t}=0,\\[5mm]
\di \mbox{whenever} \hspace{.5cm} u_t-u_{x^2t}+au_{x}(1-u_{t})=0.
\end{array} \right. \label{eq:3.4} \ee

Substituting (\ref{eq:3.3}) into (\ref{eq:3.4}), and equating the
coefficients of the various monomials in partial derivatives with
respect to $x$ and various power of $u$, we can find the
determining equations for the symmetry group of the Eq.
(\ref{eq:1.0}). Solving this equations, we get the following forms
of the coefficient functions
\be
 \xi=\frac{-x}{3}c_1+c_3, \quad \eta=c_1t+c_2, \quad  \varphi=(\frac{2t}{3}+\frac{u}{3}-\frac{2x}{3a})c_1+c_4. \label{eq:3.5} \ee
where $c_1$, $c_2$, $c_3$ and $c_4$ are arbitrary constant. Thus,
the Lie algebra $\goth g$ of infinitesimal symmetry of the Eq.
(\ref{eq:1.0}) is spanned bye the four vector fields
\be \textbf{v}_1=\partial_x,\quad \textbf{v}_2=\partial_t,\quad
\textbf{v}_3=\frac{1}{a}\partial_u,\quad
\textbf{v}_4=3t\partial_t-x\partial_x+(2t+u-\frac{2x}{a})\partial_u.
\label{eq:3.6} \ee
The commutation relations between these vector fields are given
in the Table 1. The Lie algebra $\goth g$ is solvable, because if
${\goth g}^{(1)}=\langle
\textbf{v}_i,[\textbf{v}_i,\textbf{v}_j]\rangle=[\goth g,\goth
g]$, we have ${\goth g}^{(1)}=\langle \textbf{v}_1,\cdots
\textbf{v}_4\rangle$, and ${\goth g}^(2)=[{\goth g}^{(1)},{\goth
g}^{(1)}]=\langle-\textbf{v}_1-2\textbf{v}_3,3\textbf{v}_2+2a\textbf{v}_3,\textbf{v}_3\rangle$,
so, we have a chain of ideals ${\goth g}^{(1)}\supset {\goth
g}^{(2)}\supset {0}$.

\begin{table}[h] \label{Tab:1}
  \small \caption{\textit{The commutator table}} \begin{center}
\begin{tabular}{ccccc} \hline
  $[{\mathbf v}_{i},{\mathbf v}_{j}]$ & ${\mathbf v}_1$ & ${\mathbf v}_2$ & ${\mathbf v}_3$ & ${\mathbf v}_4$ \\ \hline
  ${\mathbf v}_1$ & 0 & 0 & 0 & $-{\mathbf v}_1-2{\mathbf v}_3$ \\
  ${\mathbf v}_2$ & 0 & 0 & 0 & $3{\mathbf v}_2+2a{\mathbf v}_3$ \\
  ${\mathbf v}_3$ & 0 & 0 & 0 & ${\mathbf v}_3$\\
  ${\mathbf v}_4$ & ${\mathbf v}_1+2{\mathbf v}_3$ & $-3{\mathbf v}_2-2a{\mathbf v}_3$ & -${\mathbf v}_3$ & 0 \\\hline
\end{tabular} \end{center}
\end{table}

\medskip To obtain the group transformation which is generated by the
infinitesimal generators $\textbf{v}_i$ for $i=1,2,3,4$ we need to
solve the three systems of first order ordinary differential
equations
\be \di \frac{d\tilde{x}(s)}{ds} &=&
\xi_i(\tilde{x}(s),\tilde{t}(s),\tilde{u}(s)), \quad
\tilde{x}(0)=x, \nonumber\\
\di \frac{d\tilde{t}(s)}{ds} &=&
\eta_i(\tilde{x}(s),\tilde{t}(s),\tilde{u}(s)), \quad
\tilde{t}(0)=t, \qquad i=1,\cdots,4 \label{eq:3.7}\\
\di \frac{d\tilde{u}(s)}{ds} &=&
\varphi_i(\tilde{x}(s),\tilde{t}(s),\tilde{u}(s)), \quad
\tilde{u}(0)=u. \nonumber
 \ee
Exponentiating the infinitesimal symmetries of Eq.
(\ref{eq:1.0}), we get the one-parameter groups $G_i(s)$ generated
by $\textbf{v}_i$ for $i=1,\cdots,4$
\be
G_1:(t,x,u) & \longmapsto & (x+s,t,u),\nonumber\\
G_2:(t,x,u) & \longmapsto & (x,t+s,u),\label{eq:3.8} \\
G_3:(t,x,u) & \longmapsto & (x,t,u+s/a),\nonumber \\
G_4:(t,x,u) &\longmapsto&(x{\rm e}^{-s},t{\rm {e}}^{3s},t{\rm
{e}}^{3s}+\frac{x}{a}{\rm {e}}^{-s}+(u-t-\frac{x}{a})\rm
{{e}}^s).\nonumber \ee
Consequently,
\paragraph{Theorem 3.1}
{\it If $u=f(x,t)$ is a solution of Eq. (\ref{eq:1.0}), so are the
functions
\be G_1(s)\cdot f(x,t)&=&f(x-s,t),\nonumber\\
G_2(s)\cdot f(x,t)&=&f(x,t-s),\label{eq:3.9} \\
G_3(s)\cdot f(x,t)&=&f(x,t)+s/a,\nonumber \\
G_4(s)\cdot f(x,t)&=&{\rm {e}}^sf(x{\rm {e}}^s,t{\rm
{e}}^{-3s})+\frac{x}{a}(1-{\rm e}^{2s})+t(1-{\rm
e}^{-2s}).\nonumber \ee}
\section{Optimal system of the H-R equation }
In general, to each s-parameter subgroup $H$ of the full symmetry
group $G$ of a system of differential equations in $p>s$
independent variables, there will correspond a family of
group-invariant solutions. Since there are almost always an
infinite number of such subgroups, it is not usually feasible to
list all possible group-invariant solutions to the system. We
need an effective, systematic means of classifying these
solutions, leading to an "optimal system" of group-invariant
solutions from which every other such solution can be
derived.\cite{[3]}
\paragraph{Definition 4.1}
Let $G$ be a Lie group with Lie algebra $\goth g$. An optimal
system of $s-$parameter subgroups is a list of conjugacy
inequivalent $s-$parameter subalgebras with the property that any
other subgroup is conjugate to precisely one subgroup in the
list. Similarly, a list of $s-$parameter subalgebras forms an
optimal system if every $s-$parameter subalgebra of $\goth g$ is
equivalent to a unique member of the list under some element of
the adjoint representation: $\overline{\goth h}={\mathrm
Ad}(g({\goth h}))$.\cite{[3]}
\paragraph{Theorem 4.2}
{\it Let $H$ and $\overline{H}$ be connected s-dimensional Lie
subgroups of the Lie group $G$ with corresponding Lie subalgebras
$\goth h$ and $\overline{\goth h}$ of the Lie algebra $\goth g$ of
$G$. Then $\overline{H}$=$gHg^{-1}$ are conjugate subgroups if and
only if $\overline{\goth h}={\mathrm Ad}(g({\goth h}))$ are
conjugate subalgebras. }\cite{[3]}

\medskip By theorem (4.2), the problem of finding an optimal system of
subgroups is equivalent to that of finding an optimal system of
subalgebras. For one-dimensional subalgebras, this classification
problem is essentially the same as the problem of classifying the
orbits of the adjoint representation, since each one-dimensional
subalgebra is determined by nonzero vector in $\goth g$. This
problem is attacked by the na\"{\i}ve approach of taking a general
element ${\mathbf V}$ in $\goth g$ and subjecting it to various
adjoint transformation so as to "simplify" it as much as
possible. Thus we will deal with the construction of the optimal
system of subalgebras of $\goth g$.

\medskip To compute the adjoint representation, we use the Lie
series
\be {\mathrm Ad}(\exp(\varepsilon{\mathbf v}_i){\mathbf v}_j) =
{\mathbf v}_j-\varepsilon[{\mathbf v}_i,{\mathbf
v}_j]+\frac{\varepsilon^2}{2}[{\mathbf v}_i,[{\mathbf
v}_i,{\mathbf v}_j]]-\cdots,\label{eq:4.1} \ee
where $[{\mathbf v}_i,{\mathbf v}_j]$ is the commutator for the
Lie algebra, $\varepsilon$ is a parameter, and $i,j=1,2,3,4$. Then
we have the Table 2.
\begin{table}[h] \label{Tab:2}
\small \caption{\textit{Adjoint representation table of the
infinitesimal generators ${\mathbf v}_i$ }}
\begin{center}
\small \begin{tabular}{ccccc} \hline
  $Ad$ & ${\mathbf v}_1$ & ${\mathbf v}_2$ & ${\mathbf v}_3$ & ${\mathbf v}_4$ \\\hline
  ${\mathbf v}_1$ & ${\mathbf v}_1$ & ${\mathbf v}_2$ & ${\mathbf v}_3$ & ${\mathbf v}_4+\varepsilon({\mathbf v}_1+2{\mathbf v}_3)$ \\
  ${\mathbf v}_2$ & ${\mathbf v}_1$ & ${\mathbf v}_2$ & ${\mathbf v}_3$ & ${\mathbf v}_4-\varepsilon(3{\mathbf v}_2+2a{\mathbf v}_3)$ \\
  ${\mathbf v}_3$ & ${\mathbf v}_1$ & ${\mathbf v}_2$ & ${\mathbf v}_3$ & ${\mathbf v}_4-\varepsilon{\mathbf v}_3$ \\
  ${\mathbf v}_4$ & ${\mathbf v}_1-\varepsilon({\mathbf v}_1+2{\mathbf v}_3)$ & ${\mathbf v}_2+\varepsilon(3{\mathbf v}_2+2a{\mathbf v}_3)$ & ${\mathbf v}_3+\varepsilon{\mathbf v}_3$ & ${\mathbf v}_4$ \\\hline
\end{tabular} \end{center}
\end{table}

\paragraph{Theorem 4.3}
{ \it An optimal system of one-dimensional Lie algebras of the
H-R equation is provided by }\\
1) $\;\textbf{v}_4$, \quad 2)
$\;\alpha\textbf{v}_1+\beta\textbf{v}_2+\textbf{v}_3$, \quad
3)$\;\alpha\textbf{v}_1+\textbf{v}_2$,\quad
4)$\;\textbf{v}_1$\quad \hfill\ \mbox{ }
\paragraph{Proof:}
Consider the symmetry algebra $\goth g$ of the equation
(\ref{eq:1.0}) whose adjoint representation was determined in
table 2 and let $F^s_i:{\goth g}\to{\goth g}$ defined by
${\mathbf v}\mapsto\mathrm{Ad}(\exp(\varepsilon{\mathbf
v}_i){\mathbf v})$ is a linear map, for $i=1,\cdots,4$. The
matrices $M^\varepsilon_i$ of $F^\varepsilon_i$, $i=1,\cdots,4$,
with respect to basis $\{{\mathbf v}_1,\cdots,{\mathbf v}_4\}$ are
\begin{eqnarray*}
M^\varepsilon_1=\small \left[ \begin {array}{cccc} 1&0&0&-\varepsilon\\
0&1&0&0\\
0&0&1&-2\varepsilon\\
0&0&0&1\end {array} \right],
M^\varepsilon_2=\small \left[ \begin {array}{cccc} 1&0&0&0\\
0&1&0&3\varepsilon\\
0&0&1&2a\varepsilon\\
0&0&0&1\end {array} \right],
M^\varepsilon_3=\small \left[ \begin {array}{cccc} 1&0&0&0\\
0&1&0&0\\
0&0&1&{\varepsilon}\\
0&0&0&1
\end {array} \right],\\
M^\varepsilon_4=\small \left[ \begin {array}{cccc} {{\rm e}^{{\varepsilon}}}&0&0&0\\
0&{{\rm e}^{-3\,{\varepsilon}}}&0&0 \\
{{\rm e}^{{\varepsilon}}}-{{\rm e}^{-{\varepsilon}}}&a{
{\rm e}^{-{\varepsilon}}} \left( {{\rm e}^{-2\,{\varepsilon}}}-1 \right) &{ {\rm e}^{-{\varepsilon}}}&0\\
0&0&0&1\end {array}
 \right]
\end{eqnarray*}
Let ${\mathbf V}=\sum_{i=1}^4a_i{\mathbf v}_i$ is a nonzero
vector field in $\goth g$. We will simplify as many of the
coefficients $a_{i}$ as possible by acting these matrices on a
vector field ${\mathbf V}$ alternatively.

Suppose first that $a_{4}\neq0$, scaling ${\mathbf V}$ if
necessary we can assume that $a_{4}=1$, then we can make the
coefficients of ${\mathbf v}_1$, ${\mathbf v}_2$ and ${\mathbf
v}_3$ vanish by $M^\varepsilon_1$ and $M^\varepsilon_2$. And
${\mathbf V}$ reduced to case 1.

If $a_{4}=0$ and $a_{3}\neq0$, then we can not make vanish the
coefficients of ${\mathbf v}_1$ and ${\mathbf v}_2$ by acting any
matrices $M^\varepsilon_i$. Scaling ${\mathbf V}$ if necessary,
we can assume that $a_3=1$ and ${\mathbf V}$ reduced to case 2.

If $a_{4}=a_{3}=0$ and $a_{2}\neq0$, then we can not make vanish
the coefficient of ${\mathbf v}_1$. Scaling ${\mathbf V}$ if
necessary, we can assume that $a_2=1$ and ${\mathbf V}$ reduced
to case 3.

The remaining one-dimensional subalgebras are spanned by vectors
of the above form with $a_{4}=a_{3}=a_{2}=0$. If $a_{1}\neq0$, we
scale to make $a_{1}=1$, and  ${\mathbf V}$ reduced to case
3.\hspace{8.9cm} $\Box$
\section{Symmetry reduction of the H-R equation }
Lie-group method is applicable to both linear and non-linear
partial differential equations, which leads to similarity
variables that may be used to reduce the number of independent
variables in partial differential equations. By determining the
transformation group under which a given partial differential
equation is invariant, we can obtain information about the
invariants and symmetries of that equation.

Symmetry group method will be applied to the (\ref{eq:1.0}) to be
connected directly to some order differential equations. To do
this, a particular linear combinations of infinitesimals are
considered and their corresponding invariants are determined. The
equation (\ref{eq:1.0}) is expressed in the coordinates
$(x,t,u)$, so to reduce this equation is to search for its form
in specific coordinates. Those coordinates will be constructed by
searching for independent invariants $(y,v)$ corresponding to the
infinitesimal generator. So using the chain rule, the expression
of the equation in the new coordinate allows us to the reduced
equation. Here we will obtain some invariant solutions with
respect to symmetries. First we obtain the similarity variables
for each term of the Lie algebra $\goth g$, then we use this
method to reduced the PDE and find the invariant solutions.

We can now compute the invariants associated with the symmetry
operators, they can be obtained by integrating the characteristic
equations. For example for the operator
$\textbf{v}_2+\textbf{v}_3=\partial_t+\frac{1}{a}\partial_u$
characteristic equation is
\be \frac{dx}{0}=\frac{dt}{1}=\frac{a\,du}{1}. \label{eq:5.1} \ee
The corresponding invariants are $y=x$, $v=u-\frac{t}{a}$
therefore, a solution of our equation in this case is
$u=v(y)+\frac{t}{a}$. The derivatives of $u$ are given in terms
of $y$ and $v$ as
\be u_x=v_{y},\quad u_{x^2}=v_{yy},\quad u_{x^2t}=0,\quad
u_t=\frac{1}{a}.\label{eq:5.2}\ee
Substituting (\ref{eq:5.2}) into the Eq. (\ref{eq:1.0}), we obtain
the ordinary differential equation $(a-1)v_{y}+1/a=0$, the
solution of this equation is $v=\frac{-y}{a(a-1)}+c$.
Consequently, we obtain that
\be u=\frac{x}{a(1-a)}+\frac{t}{a}+c.\label{eq:5.3}\ee
All results are coming in the tables 3 and 4.
\begin{table}[h] \label{Tab:3}
\small \caption{\textit{Reduction of Eq. (\ref{eq:1.0}) }}
\begin{center}
\small \begin{tabular}{cccc} \hline
  operator & $y$ & $v$ & $u$ \\\hline
  ${\mathbf v}_1$ & $t$ & $u$ & $v(y)$ \\
  ${\mathbf v}_2$ & $x$ & $u$ & $v(y)$\\
  ${\mathbf v}_4$ & $xt^{1/3}$ & $(u-2x/a)t^{-1/3}+t^{2/3}$ & $v(y)t^{1/3}+2x/a-t$ \\
  ${\mathbf v}_1+a{\mathbf v}_3$ & $t$ & $u-x$ & $v(y)+x$  \\
  ${\mathbf v}_2+{\mathbf v}_3$ & $x$ & $u-t/a$ & $v(y)+t/a$ \\
  ${\mathbf v}_1+{\mathbf v}_2$ & $x-t$ & $u$ & $v(y)/a$ \\
  \hline
\end{tabular} \end{center}
\end{table}
\begin{table}[h] \label{Tab:4}
\small \caption{\textit{Reduced equations corresponding to
infinitesimal symmetries }}
\begin{center}
\small \begin{tabular}{cc} \hline
  operator & similarity reduced equations \\\hline
  ${\mathbf v}_1$ & $v_{y}=0$ \\
  ${\mathbf v}_2$ & $av_{y}=0$\\
  ${\mathbf v}_4$ & $t^{-2/3}v+(xt^{-1/3}-3)v_{y}+xt^{-2/3}v_{yy}+(at^{2/3}v_{y}+2)(6-t^{-2/3}v-t^{-1/3}v_{y})=3$ \\
  ${\mathbf v}_1+a{\mathbf v}_3$ & $(1-a)v_{y}+a=0$  \\
  ${\mathbf v}_2+{\mathbf v}_3$ & $(a-1)v_{y}+1/a=0$ \\
  ${\mathbf v}_1+{\mathbf v}_2$ & $-v_{y}+v_{yyy}+av_{y}(1+v_{y})=0$ \\
  \hline
\end{tabular} \end{center}
\end{table}

\section{Characterization of differential invariants }
Differential invariants help us to find general systems of
differential equations which admit a prescribed symmetry group.
One say, if G is a symmetry group for a system of PDEs with
functionally differential invariants, then, the system can be
rewritten in terms of differential invariants. For finding the
differential invariants of the Eq. (\ref{eq:1.0}) up to order 2,
we should solve the following systems of PDEs:
\be {\frac {\partial I}{\partial x}},\quad {\frac {\partial
I}{\partial t}},\quad \frac{1}{a}{\frac {\partial I}{\partial
u}},\quad 3t{\frac {\partial I}{\partial t}}-x{\frac {\partial
I}{\partial x}}+(2t+u-\frac{2x}{a}){\frac {\partial I}{\partial
u}}, \label{eq:6.1}\ee
where $I$ is a smooth function of $(x,t,u)$,
\be {\frac {\partial I_1}{\partial x}}, \quad {\frac {\partial
I_1}{\partial t}}, \quad \frac{1}{a}{\frac {\partial I_1}{\partial
u}},\quad 3t{\frac {\partial I_1}{\partial t}}-x{\frac {\partial
I_1}{\partial x}}+\cdots+(2-2u_{t}){\frac {\partial I_1}{\partial
u_{t}}}, \label{eq:6.2}\ee
where $I_{1}$ is a smooth function of $(x,t,u,u_{x},u_{t})$,
\be {\frac {\partial I_2}{\partial x}}, \quad {\frac {\partial
I_2}{\partial t}}, \quad \frac{1}{a}{\frac {\partial I_2}{\partial
u}},\quad 3t{\frac {\partial I_2}{\partial t}}-x{\frac {\partial
I_2}{\partial x}}+\cdots-u_{xt}{\frac {\partial I_2}{\partial
u_{xt}}}-5u_{tt}{\frac {\partial I_2}{\partial u_{tt}}},
\label{eq:6.3}\ee
where $I_{2}$ is a smooth function of
$(x,t,u,u_{x},u_{t},u_{xx},u_{xt},u_{tt})$. The solutions of PDEs
systems (\ref{eq:6.1}),(\ref{eq:6.2}) and (\ref{eq:6.3}) coming
in table 5, where * and ** are refer to ordinary invariants and
first order differential invariants respectively.
\begin{table}[h] \label{Tab:5}
\small \caption{\textit{differential invariants }}
\begin{center}
\small \begin{tabular}{|c|ccc|} \hline
  vector field & ordinary invariant  & 1st order diff. invariant & 2nd order diff. invariant \\\hline
  ${\mathbf v}_1$ & $t,u$ & $*,u_{x},u_{t}$ & $*,**,u_{xx},u_{xt},u_{tt}$ \\
  ${\mathbf v}_2$ & $x,u$ & $*,u_{x},u_{t}$ & $*,**,u_{xx},u_{xt},u_{tt}$\\
  ${\mathbf v}_3$ & $x,t$ & $*,u_{x},u_{t}$ & $*,**,u_{xx},u_{xt},u_{tt}$ \\
  ${\mathbf v}_4$ & $tx^3,(\frac{-x}{a}-t+u)x$ & $*,x^2u_{x}-\frac{1}{a},\frac{u_{t}-1}{x^2}$ & $*,**,x^3u_{xx},\frac{u_{xt}}{x},\frac{u_{tt}}{x^5}$  \\
  \hline
\end{tabular} \end{center}
\end{table}
\section{Nonclassical symmetries of the H-R equation }
In this section we would like to apply the nonclassical method to
the H-R equation. The graph of a solution
\be
u^{\alpha}=f^{\alpha}(x_{1},\cdots,x_{n}),\hspace{1.5cm}\alpha=1,\cdots,q\label{eq:7.1}\ee
to the system (\ref{eq:2.1}) defines an p-dimensional submanifold
$\Gamma_{f}\subset R^{p}\times R^{q}$ of the space of independent
and dependent variables. The solution will be invariant under the
one-parameter subgroup generated by vector (\ref{eq:2.4}) if and
only if $\Gamma_{f}$ is an invariant submanifold of this group. By
applying the well known criterion of invariance of a submanifold
under a vector field we get that (\ref{eq:7.1}) is invariant
under vector (\ref{eq:2.4}) if and only if $f$ satisfies the first
order system $E_{Q}$ of partial differential equations
\be
Q^{\alpha}(x,u,u^{(1)})=\varphi^{(\alpha)}(x,u)-\sum_{i=1}^p\xi^i(x,u)u_i^\alpha=0,\hspace{1cm}\alpha=1,\cdots,q\label{eq:7.2}\ee
known as the invariant surface conditions. The q-tuple
$Q=(Q^1,\cdots,Q^q)$ is known as the characteristic of the vector
field (\ref{eq:2.4}). In what follows, the n-th prolongation of
the invariant surface conditions (\ref{eq:7.2}) will be denoted
by $E_Q^{(n)}$, which is a n-th order system of partial
differential equations obtained by appending to (\ref{eq:7.2}) its
partial derivatives with respect to the independent variables of
orders $j\leq n-1$.

For the system (\ref{eq:2.1}), (\ref{eq:7.2}) to be compatible,
the n-th prolongation $\textrm{Pr}^{(n)}{\mathbf v}$ of the
vector field ${\mathbf v}$ must be tangent to the intersection
$E\cap E_{Q}^{(n)}$
\be \textrm{Pr}^{(n)}{\mathbf v}(\Delta_\nu)|_{E\cap
E_{Q}^{(n)}}=0, \hspace {1cm}\nu=1,\cdots,l.\label{eq:7.3}\ee
If the equations (\ref{eq:7.3}) are satisfied, then the vector
field (\ref{eq:7.3}) is called a nonclassical infinitesimal
symmetry of the system (\ref{eq:2.1}). The relations
(\ref{eq:7.3}) are generalizations of the relations
(\ref{eq:2.5}) for the vector fields of the infinitesimal
classical symmetries. A similar procedure is applicable to the
case of the nonclassical infinitesimal symmetries with an evident
difference that in general one has fewer determining equations
than in the classical case. Therefore, we expect that nonclassical
symmetries are much more numerous than classical ones, since any
classical symmetry is clearly a nonclassical one. The important
feature of determining equations for nonclassical symmetries is
that they are nonlinear, this implies that the space of
nonclassical symmetries does not, in general, form a vector
space. For more theoretical background see \cite{[20],[5]}.

Consider the system $E$ of second order equations
\be u_t-v_t+au_x(1-u_t)=0, \hspace{1cm} u_{xx}-v=0
\label{eq:7.4}\ee
obtained from the H-R equation. If we assume that the coefficient
of $\partial_t$ of the vector field (\ref{eq:2.4}) does not
identically equal zero, then for the vector field
\be {\mathbf
v}=\xi(x,t,u,v)\partial_x+\partial_t+\varphi(x,t,u,v)\partial_u+\psi(x,t,u,v)\partial_v
\label{eq:7.5}\ee
the invariant surface conditions are
\be u_t+\xi u_x=\varphi, \hspace{1cm} v_t+\xi v_x=\psi
\label{eq:7.6}\ee
The equations (\ref{eq:7.3}) take the forms
\be \begin{array}{lcl} (2\,(u_{t}-1/2)u_x^{2}a-u_t(u_{x}-v_x))\xi_u-u_{t}\psi_u-\psi_t \\
+((1-2\,u_t)au_x+u_t)\varphi_u+(u_x(1-u_t)a-u_t)\psi_v \\
+(-u_{x}^3( u_{t}-1)a^2+2\,(u_{t}-1/2)u_{x}^2a-u_t( u_{x}-v_{x}))\xi_v \\
+(u_{x}^2(u_{t}-1)a^2+((1-2\,u_{t})u_{x}-v_{x}(u_{t}-1))a+u_{t})\varphi_{v} \\
+(1-au_x)\varphi_t-a(u_t-1)\varphi_x+u_x(u_t-1)a\xi_x+(u_{x}^2a-u_x+v_x)\xi_t=0,
\end{array}\label{eq:7.7}\ee
and
\be \begin{array}{lcl} -\psi-u_x\xi_{xx}-2\,u_x^2\xi_{xu}
+u_x^2\varphi_{uu}-u_x^3\xi_{uu}+u_{xx}\varphi_u
-2\,u_{xx}\xi_{x}+v_x^2\varphi_{vv}\\
+v_{xx}\varphi_{v}+2\,u_x\varphi_{xu}+2\,v_x\varphi_{xv}-2\,u_xv_x\xi_{xv}
+2\,u_xv_x\varphi_{uv}-2\,u_x^2v_x\xi_{uv} \\
-3\,u_{xx}u_x\xi_{u}-2\,u_{xx}v_x\xi_{v}-u_xv_x^2\xi_{vv}-v_{xx}u_x\xi_v+\varphi_{xx}=0.
\end{array}\label{eq:7.8}\ee
After inserting $\psi$ and its derivatives, as determined by the
equation (\ref{eq:7.8}), in to (\ref{eq:7.7}) and substituting
$v_x=u_{xxx}, v_{xx}=u_{xxxx}$, and equating the coefficients of
the various monomials in partial derivatives with respect to $x$
and various power of $u$, we can find the determining equations.
Solving this equations, we get four Algebraic equations equal to
zero. This means that no supplementary symmetries, of
non-classical type, are specific for our model.

Now assume that the coefficient of $\partial_t$ in (\ref{eq:7.5})
equals zero and try to find the infinitesimal nonclassical
symmetries of the form
\be {\mathbf
v}=\partial_x+\varphi(x,t,u,v)\partial_u+\psi(x,t,u,v)\partial_v
\label{eq:7.9}\ee
for which the invariant surface conditions are the following ones
\be u_x=\varphi, \hspace{1cm} v_x=\psi \label{eq:7.10}\ee
Relations (\ref{eq:7.3}) lead to the system of equations for the
functions $\varphi$ and $\psi$
\be \begin{array}{lcl}
(u_x^2(u_t-1)a^2+((-v_x-2\,u_x)u_t+v_x+u_x)a+ u_t)\varphi_v\\+
((u_x-2\,u_xu_t)a+u_t)\varphi_u+(u_x(u_t-1)a-u_t)\psi_v \\
+(1-au_x)\varphi_t-a(u_t-1)\varphi_x-\psi_t -u_t\psi_u=0,
\end{array}\label{eq:7.11}\ee
and
\be \begin{array}{lcl} -\psi+\varphi_{xx}+2\,u_x\varphi_{xu}+2\,
\varphi_{xv}v_x+u_x^2\varphi_{uu}+2\,u_xv_x\varphi_{uv}\\
+u_{xx}\varphi_u+v_x^2\varphi_{vv}+ v_{xx}\varphi_{v}=0.
\end{array}\label{eq:7.12}\ee
Similar the previous case, we can find determining equations.
Solving this equations, we get the following form of the
coefficient functions
\be \varphi=c_1u-\frac{c_1x}{a}+t(1-2u_t)c_1+c_2
\label{eq:7.13}\ee
where $c_1$ and $c_2$ are arbitrary constant. So the system
(\ref{eq:7.4}) admits the classical symmetry ${\mathbf v_1}$ and
nonclassical symmetry ${\mathbf v_2}$:
\be {\mathbf v_1}=\partial_x+\partial_u, \hspace{1cm} {\mathbf
v_2}=\partial_x-(\frac{x}{a}-u-t+2\,tu_t)\partial_u
\label{eq:7.14}\ee
\section{Conservation laws of the H-R equation  }
Many methods for dealing with the conservation laws are derived,
such as the method based on the Noether's theorem, the multiplier
method, by the relationship between the conserved vector of a PDE
and the Lie-Bäcklund symmetry generators of the PDE, the direct
method, etc.\cite{[3],[16],[17],[18]}.

Now, we derive the conservation laws from the multiplier method.
\paragraph{Definition 8.1}
A local conservation law of the PDE system (\ref{eq:2.1}) is a
divergence expression
\be D_i\Phi^i[u]=D_1\Phi^1[u]+\cdots+D_n\Phi^n[u]=0
\label{eq:8.1}\ee
holding for all solutions of the system (\ref{eq:2.1}). In
(\ref{eq:8.1}),
$\Phi^i[u]=\Phi^i(x,u,\partial_u,\cdots,\partial^r_u)$,
$i=1,\cdots,n$, are called fluxes of the conservation law, and
the highest-order derivative $(r)$ present in the fluxes
$\Phi^i[u]$ is called the order of a conservation law. \cite{[17]}

\medskip In particular, a set of multipliers
$\{\Lambda_\nu[U]\}^l_{\nu=1}=\{\Lambda_\nu(x,U,\partial_U,\cdots,\partial^r_U)\}^l_{\nu=1}$
yields a divergence expression for the system
$\Delta_\nu(x,u^{(n)})$ such that if the identity
\be \Lambda_\nu[U]\Delta_\nu[U] \equiv D_i\Phi^i[U]
\label{eq:8.2}\ee
holds identically for arbitrary functions $U(x)$. Then on the
solutions $U(x)=u(x)$ of the system (\ref{eq:2.1}), if
$\Lambda_\nu[U]$ is non-singular, one has local conservation law
$\Lambda_\nu[u]\Delta_\nu[u]=D_i\Phi^i[u]=0$.
\paragraph{Definition 8.2}
The Euler operator with respect to $U^j$ is the operator defined
by
\be
E_{U^j}=\frac{\partial}{\partial{U^j}}-D_i\frac{\partial}{\partial{U^j}}+\cdots+(-1)^sD_{i_1}\cdots
D_{i_s}\frac{\partial}{\partial{U^j_{i_1\cdots i_s}}}+\cdots
\label{eq:8.3}\ee
for $j=1,\cdots,q$. \cite{[17]}

\paragraph{Theorem 8.3}
The equations
$E_{U^j}F(x,U,\partial_U,\cdots,\partial^s_U)\equiv0$,
$j=1,\cdots,q$ hold for arbitrary $U(x)$ if and only if
$F(x,U,\partial_U,\cdots,\partial^s_U)\equiv
D_i\Psi^i(x,U,\partial_U,\cdots,\partial^{s-1}_U)$ holds for some
functions $\Psi^i(x,U,\partial_U,\cdots,\partial^{s-1}_U)$,\,
$i=1,\cdots q$. \cite{[17]}\hspace{2.5cm} $\Box$
\paragraph{Theorem 8.4}
A set of non-singular local multipliers
$\{\Lambda_\nu(x,U,\partial_U,\cdots,\partial^r_U)\}^l_{\nu=1}$
yields a local conservation law for the system
$\Delta_\nu(x,u^{(n)})$ if and only if the set of identities
\be
E_{U^j}(\Lambda_\nu(x,U,\partial_U,\cdots,\partial^r_U)\Delta_\nu(x,u^{(n)}))\equiv
0, \,j=1,\cdots q,\label{eq:8.4}\ee
holds for arbitrary functions $U(x)$. \cite{[17]}\hspace{6.2cm}
$\Box$

\medskip The set of equations (\ref{eq:8.4}) yields the set of linear
determining equations to find all sets of local conservation law
multipliers of the system (\ref{eq:2.1}). Now, we seek all local
conservation law multipliers of the form
$\Lambda=\xi(x,t,u,u_x,u_t,u_{xx},u_{xt},u_{tt})$ of the equation
(\ref{eq:1.0}). The determining equations (\ref{eq:8.4}) become
\be
E_U[\xi(x,t,U,U_x,U_t,U_{xx},U_{xt},U_{tt})(U_t-U_{x^2t}+aU_{x}(1-U_{t}))]\equiv
0,\label{eq:8.5}\ee
where $U(x,t)$ are arbitrary function. Equation (\ref{eq:8.5})
split with respect to third order derivatives of $U$ to yield the
determining PDE system whose solutions are the sets of local
multipliers of all nontrivial local conservation laws of second
order of H-R equation.

The solution of the determining system (\ref{eq:8.5}) given by
\be c_1\,U_{xx}+\frac{1}{2}\,c_2(
2\,tU_{tt}+U_{t}-1)+c_3u_{tt},\label{eq:8.6}\ee
where $c_1$, $c_2$ and $c_3$ are arbitrary constant. So local
multipliers given by
\be 1)\; \xi=U_{xx}, \hspace{1cm} 2)\; \xi=U_{tt}, \hspace{1cm}
3)\; \xi=tU_{tt}+\frac{1}{2}U_{t}-\frac{1}{2},\label{eq:8.7}\ee
Each of the local multipliers $\xi$ determines a nontrivial
two-order local conservation law $D_t\Psi+D_x\Phi=0$ with the
characteristic form
\be D_t\Psi+D_x\Phi \equiv
\xi(U_t-U_{x^2t}+aU_{x}(1-U_{t})),\label{eq:8.8}\ee
To calculate the conserved quantities $\Psi$ and $\Phi$, we need
to invert the total divergence operator. This requires the
integration (by parts) of an expression in multi-dimensions
involving arbitrary functions and its derivatives, which is a
difficult and cumbersome task. The homotopy operator \cite{[19]}
is a powerful algorithmic tool (explicit formula) that originates
from homological algebra and variational bi-complexes.

\paragraph{Definition 8.5}
The 2-dimensional homotopy operator is a vector operator with two
components,
$\Big({\textit{H}}^{(x)}_{u(x,t)}\textit{f},{\textit{H}}^{(t)}_{u(x,t)}\textit{f}\Big)$,
where
\be
{\textit{H}}^{(x)}_{u(x,t)}\textit{f}=\int^1_0\Big(\sum^q_{j=1}\textit{I}^{(x)}_{u^j}\textit{f}\Big)[\lambda
u]\frac{d\lambda}{\lambda}\;\;\;\hspace{.2cm} and\;\;\;
\hspace{.2cm}
{\textit{H}}^{(t)}_{u(x,t)}\textit{f}=\int^1_0\Big(\sum^q_{j=1}\textit{I}^{(t)}_{u^j}\textit{f}\Big)[\lambda
u]\frac{d\lambda}{\lambda}.\label{eq:8.9}\ee
The x-integrand, $\textit{I}^{(x)}_{u^j_{(x,t)}}\textit{f}$, is
given by
\be
\textit{I}^{(x)}_{u^j}\textit{f}=\sum^{M^j_1}_{k_1=1}\sum^{M^j_2}_{k_2=0}\Big(\sum^{k_1-1}_{i_1=0}\sum^{k_2}_{i_2=0}
B^{(x)}u^j_{x^{i_1}t^{i_2}}(-D_x)^{k_1-i_1-1}(-D_t)^{k_2-i_2}\Big)\frac{\partial\textit{f}}{\partial
u^j_{x^{k_1}t^{k_2}}},\label{eq:8.10}\ee
where $M^j_1$, $M^j_2$ are the order of $\textit{f}$ in $u$ to
$x$ and $t$ respectively and combinatorial coefficient
\be
B^{(x)}=B(i_1,i_2,k_1,k_2)=\footnotesize{\frac{\left(\!\!\!\!\begin{array}{c}i_1+i_2\\i_1
\end{array}\!\!\!\!\right)\left(\!\!\!\!\begin{array}{c}k_1+k_2-i_1-i_2-1\\k_1-i_1-1
\end{array}\!\!\!\!\right)}{\left(\!\!\!\!\begin{array}{c}k_1+k_2\\k_1
\end{array}\!\!\!\!\right)}}.
\label{eq:8.11}\ee
Similarly, the t-integrand,
$\textit{I}^{(t)}_{u^j_{(x,t)}}\textit{f}$, is defined as
\be
\textit{I}^{(t)}_{u^j}\textit{f}=\sum^{M^j_1}_{k_1=0}\sum^{M^j_2}_{k_2=1}\Big(\sum^{k_1}_{i_1=0}\sum^{k_2-1}_{i_2=0}
B^{(t)}u^j_{x^{i_1}t^{i_2}}(-D_x)^{k_1-i_1}(-D_t)^{k_2-i_2-1}\Big)\frac{\partial\textit{f}}{\partial
u^j_{x^{k_1}t^{k_2}}},\label{eq:8.12}\ee
where $B^{(t)}(i_2,i_1,k_2,k_1)$.

\medskip For instance we apply homotopy operator to find conserved
quantities $\Psi$ and $\Phi$ which yield of multiplier
$\xi=u_{tt}$. We have
\be
\textit{f}=u_{tt}(u_t-u_{x^2t}+au_{x}(1-u_{t})),\label{eq:8.13}\ee
the integrands (\ref{eq:8.10}) and (\ref{eq:8.12}) are
\be \begin{array}{lcl}
\textit{I}^{(x)}_{u^j}\textit{f}=auu_{t^2}-auu_{t}u_{t^2}-\frac{2}{3}uu_{xt^3}+\frac{1}{3}u_{t}u_{xt^2}+\frac{1}{3}u_{x}u_{t^3}-\frac{2}{3}u_{xt}u_{t^2},\\
\\
\textit{I}^{(t)}_{u^j}\textit{f}=u^2_t-u_{t}u_{x^2t}+au_{x}u_{t}-au_{x}u^2_{t}+\frac{2}{3}uu_{x^2t^2}-auu_{xt}+auu_{xt}u_{t}\\
\hspace{1cm}+\frac{1}{3}u_{x}u_{xt^2}-\frac{1}{3}u_{x^2}u_{t^2},\end{array}\label{eq:8.14}\ee
apply (\ref{eq:8.9}) to the integrands (\ref{eq:8.14}), therefore
\be \begin{array}{lcl} {\textit{H}}^{(x)}_{u(x,t)}\textit{f}=\frac{1}{2}auu_{t^2}-\frac{1}{3}auu_{t}u_{t^2}-\frac{1}{3}uu_{xt^3}+\frac{1}{6}u_{t}u_{xt^2}
+\frac{1}{6}u_{x}u_{t^3}-\frac{1}{3}u_{xt}u_{t^2},\\
\\
{\textit{H}}^{(t)}_{u(x,t)}\textit{f}=\frac{1}{2}u^2_t-\frac{1}{2}u_{t}u_{x^2t}+\frac{1}{2}au_{x}u_{t}-\frac{1}{3}au_{x}u^2_{t}+\frac{1}{3}uu_{x^2t^2}-\frac{1}{2}auu_{xt}\\
\hspace{1.5cm}+\frac{1}{3}auu_{xt}u_{t}+\frac{1}{6}u_{x}u_{xt^2}-\frac{1}{6}u_{x^2}u_{t^2},
\end{array}\label{eq:8.15}\ee
so, we have the first conservation low of the H-R equation respect
to multiplier $\xi=u_{tt}$
\be \begin{array}{lcl}
D_x\big(\frac{1}{2}auu_{t^2}-\frac{1}{3}auu_{t}u_{t^2}-\frac{1}{3}uu_{xt^3}
+\frac{1}{6}u_{t}u_{xt^2}+\frac{1}{6}u_{x}u_{t^3}-\frac{1}{3}u_{xt}u_{t^2}\big)\\
\\
+D_t\big(\frac{1}{2}u^2_t-\frac{1}{2}u_{t}u_{x^2t}+\frac{1}{2}au_{x}u_{t}-\frac{1}{3}au_{x}u^2_{t}+\frac{1}{3}uu_{x^2t^2}-\frac{1}{2}auu_{xt}\\
\hspace{.8cm}+\frac{1}{3}auu_{xt}u_{t}+\frac{1}{6}u_{x}u_{xt^2}-\frac{1}{6}u_{x^2}u_{t^2}\big)=0.
\end{array}\label{eq:8.16}\ee
Similarly, conservation law respect to multiplier $\xi=u_{xx}$ is
\be \begin{array}{lcl}
D_x\big(\frac{1}{2}u_{x}u_{t}+\frac{1}{2}au^2_{x}-\frac{1}{3}au^2_{x}u_{t}-\frac{1}{2}uu_{xt}+\frac{1}{3}auu_{x}u_{xt})\\
+D_t(\frac{1}{2}uu_{xx}-\frac{1}{3}auu_{x}u_{t^2}-\frac{1}{2}u_{xx}\big)=0.
\end{array}\label{eq:8.17}\ee
and conservation law respect to multiplier
$\xi=tu_{tt}+\frac{1}{2}u_{t}-\frac{1}{2}$ is
\be \begin{array}{lcl}
D_x\big(-\frac{1}{2}au+\frac{1}{2}auu_{t}+\frac{1}{2}atuu_{t^2}-\frac{1}{3}atuu_{t}u_{t^2}-\frac{1}{6}auu^2_{t}-\frac{1}{2}uu_{xt^2}\\
\hspace{.7cm}-\frac{1}{3}tuu_{xt^3}+\frac{1}{6}tu_{t}u_{xt^2}-\frac{1}{12}u_{t}u_{xt}+\frac{1}{4}u_{x}u_{t^2}+\frac{1}{6}tu_{x}u_{t^3}+\frac{1}{3}u_{xt}\\
\hspace{.7cm}-\frac{1}{3}tu_{xt}u_{t^2}\big)=0.\\
+D_t\big(-\frac{1}{2}u+\frac{1}{6}uu_{x^2t}+\frac{1}{3}tuu_{x^2t^2}+\frac{1}{6}tu_xu_{xt^2}+\frac{1}{12}u_xu_{xt}+\frac{1}{6}u_{x^2}\\
\hspace{.9cm}-\frac{1}{12}u_{x^2}u_{t}-\frac{1}{2}atuu_{xt}+\frac{1}{3}atuu_{t}u_{xt}+\frac{1}{2}tu^2_{t}-\frac{1}{2}tu_tu_{x^2t}\\
\hspace{.9cm}-\frac{1}{6}tu_{x^2}u_{t^2}+\frac{1}{2}atu_xu_{t}-\frac{1}{3}atu_{x}u^2_{t}\big)=0.
\end{array}\label{eq:8.18}\ee

\end{document}